\magnification=\magstep1
\input epsf

\def\mapright#1{\smash{
	\mathop{\longrightarrow}\limits^{#1}}}

~~\vskip .5in

\centerline{\bf Decorated Teichm\"uller Theory}

\centerline{\bf  of}

\centerline{\bf Bordered Surfaces}

\vskip .5in

\centerline {\bf R. C. Penner}

\vskip .2in

\centerline{Departments of Mathematics and Physics/Astronomy}

\centerline{University of Southern California}

\centerline {Los Angeles, CA 90089}

\vskip .4in

\leftskip .8in\rightskip .8in

\noindent {\bf Abstract}~~
This paper extends the decorated Teichm\"uller theory developed before for punctured surfaces to the setting of ``bordered'' surfaces,
i.e., surfaces with boundary, and there is non-trivial new structure discovered.  The main new result identifies the arc complex of a
bordered surface up to proper homotopy equivalence with a certain quotient of the moduli space, namely, the quotient by
the natural action of the positive reals by homothety on the hyperbolic lengths of geodesic boundary components.  One tool in the proof
is a homeomorphism between two versions of a ``decorated'' moduli space for bordered surfaces.  The explicit homeomorphism relies upon
points equidistant to suitable triples of horocycles.  

\leftskip=0ex\rightskip=0ex

~~\vskip .4in

\centerline{ \bf Introduction}

\vskip .2in

\noindent  Complexes of arc families in surfaces arise in several
related contexts in mathematics.  Poincar\'e dual to an arc family in a surface is a graph embedded in the surface, so complexes of arc
families are also manifest as suitable spaces of graphs.  Such arc or graph complexes arise in many related mathematical contexts
in the works of Culler-Vogtmann, Harer, Kontsevich, the author, Strebel, and others.  Up to this point, these graphical techniques for
Riemann surfaces have been utilized in the setting of punctured surfaces without boundary.  The main purpose of this paper is to present
the analogous theory for surfaces with boundary, or so-called ``bordered'' surfaces, so in effect, we present the relative version of the
established theory for Riemann surfaces.  

\vskip .1in

\noindent In the punctured case, the quotient of an open dense subspace of the arc complex by the mapping class group is homeomorphic to
moduli space.  Our main results in this paper analogously extend to the bordered case and give a proper homotopy
equivalence between an open dense suspace  of the arc complex of a bordered surface and an ${\bf R}_+$-quotient of the moduli space of
the bordered surface (to be defined later).  An essential difference between the punctured and bordered cases of arc complexes is that in
the latter case the arcs in an arc family come in a natural linear ordering; in effect, this kills all finite isotropy in the action of
the mapping class group. 

\vskip .1in

\noindent The arc complex itself thus forms a combinatorial compactification of the moduli space in the punctured case and of its
quotient by
${\bf R}_+$ in the bordered case.  These compactifications were studied in [8], where it was conjectured that this compactification is an
orbifold in the punctured case and a sphere in the bordered case.  (In the established theory for punctured surfaces, there are other
known sphericity  results, but these are on the level of 
families of arc in a fixed surface, not its quotient by the mapping class group as in the current case of bordered surfaces.)
This conjecture was proved for the case of multiply punctured spheres in [8], as follows from the sphericity conjecture for the
relatively simple case of multiply punctured disks in the bordered case.

\vskip .1in

\noindent The proper homotopy equivalence mentioned before is produced from an identification of two versions of a ``decorated'' moduli
space; one version is based upon surfaces with horocycles about cusps in analogy to the punctured case, and the other is based upon
pairs of distinct labeled points in the geodesic boundary of a hyperbolic surface.  The identification of the two versions depends upon a
construction based upon equidistant points to triples of horocycles, which is the heart of this paper. 

\vskip .1in

\noindent Let $F=F_{g,r}^s$ denote a smooth surface of genus $g\geq 0$ with $r\geq 0$ labeled boundary components and $s\geq 0$ labeled
punctures, where 
$6g-7+4r+2s\geq 0$ (so we exclude only the surfaces $F_{1,0}^0$ and $F_{0,1}^1$).  The {\it pure mapping class group} $PMC=PMC(F)$
of
$F$ is the group of all isotopy classes of orientation-preserving homeomorphisms of $F$ pointwise fixing each boundary component
and each puncture, where the isotopy is likewise required to pointwise fix these sets. 

\vskip .1in

\noindent  

\vskip .1in

\noindent To begin with the case $r=0$, let us define the classical 
{\it Teichm\"uller space} $T_{g}^s$ of $F$ to be the
space of all complete finite-area metrics of constant Gauss
curvature ${-1}$ (so-called ``hyperbolic metrics'') on $F$ 
modulo push-forward by diffeomorphisms fixing each puncture which are isotopic to
the identity relative to the punctures.   As is well
known, $T_{g}^s$ is homeomorphic to an open ball of real dimension
$6g-6+2s$, $PMC$ acts on $T_{g}^s$ by push-forward of metric under a
representative diffeomorphism, and the
quotient $M_{g}^s=T_{g}^s/PMC$ is {\it Riemann's
moduli space} of $F$.  There is a trivial ${\bf R}_+^s$-bundle $\widetilde {\cal T}_{g}^s\to T_{g}^s$, where
the fiber over a point is the collection of all $s$-tuples of horocycles in $F$, one horocycle about each puncture, and this
leads to the corresponding trivial bundle $\widetilde{\cal T}_g^s/PMC=\widetilde {M}_{g}^s\to M_{g}^s$; the total spaces of
these bundles are respectively called the {\it decorated Teichm\"uller} and {\it decorated moduli spaces}, which are studied in
[6-10].

\vskip .1in

\noindent Turning to the case $r\neq 0$ of bordered surfaces, there are two essentially different geometric treatments of a
distinguished point
$\xi$ in a boundary component
$C$ of a surface, which lead to two different versions for the corresponding (decorated) moduli space for bordered surfaces.  In the
first treatment, we take a hyperbolic metric so that $C\ni \xi$ is a geodesic as illustrated in Figure 1a, and this leads to
one version of the moduli space of a bordered surface as we shall see (in $\S$1).  In the second treatment, we
remove the distinguished point
$\xi$ from
$C$ and take a hyperbolic metric so that $C-\{ \xi\}$ is totally geodesic as illustrated in Figure 1b (with a more hyperbolically
realistic depiction given in Figure~2), and this leads to alternate versions of moduli space
as we shall also see (in
$\S$1).  We enrich the
structure in the first treatment by ``decorating'' with a second point $p\in C$, where $p\neq \xi$, and in the second treatment by
specifying a horocyclic segment $h$ centered at
$\xi$, as are also illustrated in Figure 1, and in this way we define two versions for the decorated moduli space of a bordered
surface. After a small retraction, we prove that these two versions of decorated moduli space are homeomorphic.

\hskip 1.in
{{{\epsfysize2in\epsffile{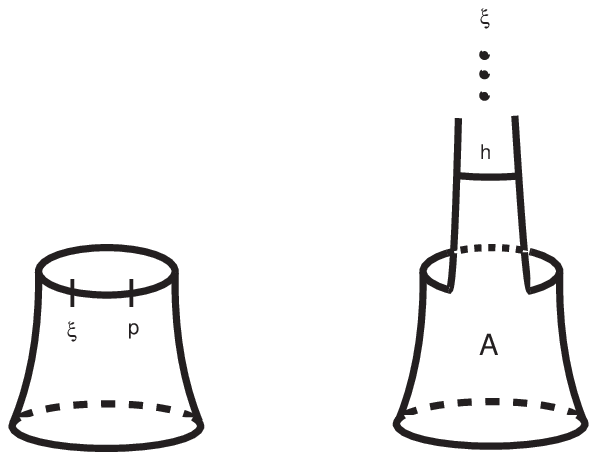}}}}

\vskip .2in

\hskip .7in {\bf 1a} Geodesic boundary. \hskip .1in{\bf 1b} Totally geodesic boundary.

\centerline{{\bf Figure 1}~~{Different treatments for points in the boundary of a surface.}}

\vskip .4in

\noindent The proof that the two versions are homeomorphic depends upon new applications (in $\S$5) of the constructions and calculations
in [7], where we studied points in hyperbolic space which are equidistant to tuples of horocycles.  The proof further involves the
extension of the decorated Teichm\"uller theory for punctured unbordered surfaces, which was developed in [6], to the setting
of bordered surfaces; we shall find that many of the arguments in [6] extend painlessly to the current case.  
There is no doubt that the serious reader must
consult the paper [6] (and perhaps [7] as well) for complete details of some of the arguments given here, but we have endeavored to keep
this note logically self-contained and complete by recalling here enough of the relevant material.

\vskip .1in

\noindent Another purpose of this paper is to present the decorated Teichm\"uller theory for bordered surfaces (in $\S$3): We
shall give both global ``lambda length'' coordinates on the decorated Teichm\"uller space as well as a $PMC$-invariant cell
decomposition of it based on a ``convex hull construction''.  In the literature for the case of bordered surfaces, Kojima has
previously described a related convex hull construction in [3], and Chekhov and Fock [1] have previously given analogous
global coordinates on Teichm\"uller space.  We shall also in an extended side-remark comment on further extensions for bordered
surfaces of the decorated Teichm\"uller theory which are interesting but are not needed here.

\vskip .1in

\noindent The main purpose of this paper is to understand the arc complex of a bordered surface (defined in $\S$6).  There is
a distinguished subspace of the arc complex which  corresponds to arc families which ``fill'' $F$ in a precise sense.  We
shall prove (in $\S$6) that this subspace is proper homotopy equivalent to the quotient of the moduli space of the bordered surface 
(in the version with geodesic boundary) by the natural ${\bf R}_+$-action by homothety on the hyperbolic lengths of the boundary
geodesics.

\vskip .1in

\noindent This paper is organized as follows.  $\S$1 defines two versions for decorated moduli space, and $\S$2
contains relevant definitions regarding arcs in bordered surfaces.  The extension of the decorated
Teichm\"uller theory to bordered surfaces is described in $\S$3, where we both recall certain arguments from  [6] for completeness
and sketch extensions of other arguments from [6] with technical details.  $\S$4 is dedicated to the study of points equidistant to
triples of horocycles, with elementary calculations from [7] simply recalled and not re-proved here, as well as a deformation
retraction of one version of decorated moduli space which is required in the sequel.  The previous material is applied in $\S$5 to give
the real-analytic homeomorphism between the two versions of decorated moduli space.  Circle actions are studied
in $\S$6 and the arc complex is defined; the proper homotopy equivalence between the ``filling'' subspace of the arc complex and the
quotient of moduli space by the homothetic ${\bf R}_+$-action is also presented in $\S$6.

\vskip .1in

\noindent To close this Introduction before turning to bordered surfaces in the sequel, we shall contrast some of the main
constructions and results in [6-10] (``in the hyperbolic setting'') for unbordered punctured surfaces
with results and constructions (``in the conformal setting'') using quadratic differentials.  In the conformal setting, Riemann's moduli
space of $F$ is regarded as the space of all equivalence classes of conformal structures on $F$. 

\vskip .1in

\noindent We begin in the conformal setting and describe the $PMC(F)$-invariant cell decomposition of $\widetilde{\cal T}_g^s$,
which is due to Harer-Mumford-Thurston [2], and relies on the Jenkins-Strebel theory [12].  The formulation of the combinatorics given
here relies on graphs with extra structure as in [10], [5].

\vskip .1in

\noindent  A ``fatgraph'' or ``ribbon graph'' $G$ is a graph with vertices at least tri-valent together with a cyclic ordering on the
half-edges about each vertex.
$G$ may be ``fattened'' to a bordered surface in the 
following way: begin with disjoint planar neighborhoods of the vertices of $G$, where the cyclic ordering agrees with that induced by
the orientation of the plane;  glue orientation-preserving bands, one band for each edge of $G$, in the natural way to these
neighborhoods.  This produces a topological surface $F_G\supseteq G$, where $G$ is a spine of $F_G$.
The complement $F_G-G$ is a collection of topological
annuli
$A_i\subseteq F_G$, $i=1,\ldots ,s\geq 1$; each annulus $A_i$ has one boundary component $\partial _i'$ lying in $G$ and the other
$\partial _i$ in
$\partial F_G$.

\vskip .1in

\noindent Let $E(G)$ denote the set of edges of $G$.
A ``metric'' on a fatgraph is a function $w\in ({\bf R}_+\cup\{ 0\} )^{E(G)}$.   
The curve $\partial _i'$ a closed edge-path on
$G$, and we define the ``length'' of $A_i$ to be $\ell _i (w)=\sum w(e)$, where the sum is over $e\in\partial _i'$ counted
with multiplicity.  A metric $w$ is ``positive'' if every length $\ell _i (w)$ is positive; let $\sigma (G)$ denote the space of all
positive metrics on $G$.

\vskip .1in

\noindent Given a positive metric $w\in \sigma (G)$ with corresponding lengths $\ell _i>0$, for $i=1,\ldots ,s$, we may construct a
metric surface homeomorphic to $F_g^s$ as follows, where $2-2g-s$ is the Euler characteristic of $G$.  Give each $A_i=A_i(w)$ the
structure of a flat cylinder with circumference
$\ell _i$ and height unity, so $A_i$ has modulus $\ell _i$; isometrically identify these cylinders in the natural way along common
edges in
$\cup\{
\partial _i'\} _1^s$ as dictated by the metric and fatgraph. This produces a metric structure on $F_G$, where the boundary component
$\partial _i$ of
$F_G$ is a standard circle of circumference $\ell _i$; glue to each $\partial _i$ a standard flat disk of
circumference $\ell _i$ with puncture $*_i$, where the boundary of the disk is concentric with $*_i$, for each $i=1,\ldots ,s$.  This
produces a conformal structure on a surface which we may identify with $F_g^s$. 

\vskip .1in

\noindent An analytic fact is:

\vskip .2in

\noindent {\bf Theorem A}~ [12;$\S$23.5]~~\it Given a positive metric $w$ on a fatgraph $G$ with lengths $\ell _i=\ell _i (w)$, for
$i=1,\ldots ,s$, there is a unique meromorphic quadratic differential $q$ on $F_g^s$ so that for each $i=1,\ldots ,s$: 

\vskip .1in

\leftskip .2in

\noindent the non-critical horizontal
trajectories of $q$ in $F$ foliate $A_i(w)\subseteq F_g^s$ by curves homotopic to the cores;

\vskip .1in

\noindent the residue of $\sqrt{q}$ at $*_i$ is $\ell _i$.\rm

\leftskip=0ex\rm

\vskip .2in

\noindent Let $\mu _q$ denote the conformal structure on $F_g^s$ determined by $q$.  We may think  
of Theorem~A abstractly as a mapping $$(G,w )\mapsto \mu _q \times (\ell _i )_1^s\in{ T}_g^s\times {\bf R}_+^s,$$
where the effective construction of $\mu _q \times (\ell _i )_1^s$ was described before.

\vskip .1in

\noindent Now, fix a surface $F=F_g^s$ and consider the collection $C_g^s$ of all homotopy classes of inclusions $G\subseteq
F$, where $G$ is a strong deformation retraction of $F$. If $G\subseteq F$ and $w\in\sigma (G)$, then we may produce
another $G_w\subseteq F$ by contracting each edge $e\in E(G)$ with $w(e)=0$ to produce $G_w$. Identifying $E(G_w)$ with 
$\{ e\in E(G): w(e)\neq 0\}$ in the natural way, we may also induce
$w'\in \sigma (G_w)$ by requiring that $w'(e)=w(e)$, for any $e\in E(G)$ with $w(e)\neq 0$.

\vskip .1in

\noindent Define 
$$U_g^s=\biggl [~\coprod _{(G\subseteq F)\in C_g^s}{\sigma (G)}~\biggr ]/\sim,$$
where $\coprod$ denotes disjoint union, and $(G^1,w ^1)\sim (G^2,w ^2)$ if and only if $G^1_{w_1}\subseteq F$ agrees with
$G^2_{w_2}\subseteq F$ as members of
$C_g^s$ and $w_1'\in\sigma (G^1)$ agrees with $w_2'\in\sigma (G^2)$. 
$PMC(F)$ acts on $U_g^s$ in the natural way induced by $(\phi:G\to F)\mapsto (f\circ\phi :G\to F)$ if $f:F\to F$ is a homeomorphism.
Define the ``fatgraph complex''
$G_g^s=U_g^s/PMC$ and let $[G,w]\in G_g^s$ denote the class of $(G,w)\in U_g^s$.

\vskip .2in

\noindent{\bf Theorem B} [12;$\S$25.6] [11]~\it ~The mapping $(G,w )\mapsto \mu _q \times (\ell _i)_1^s$ induces real-analytic
homeomorphisms $
U_g^s\to { T}_g^s\times {\bf R}_+^s$ and 
$G_g^s\to M_g^s\times {\bf R}_+^s$.\rm

\vskip .2in

\noindent  Theorem~A says that the mapping $G_g^s\to M_g^s\times {\bf R}_+^s$ is well-defined and one-to-one, while the
further analytic content of Theorem B is that this mapping is moreover onto.  The inverse map $T_g^s\times{\bf R}_+^s\to U_g^s$ is
transcendental and highly non-computable.  Theorem B gives a $PMC$-invariant cell decomposition of $T_g^s\times{\bf R}_+^s$ induced by
the cell structure of $G_g^s$.

\vskip .1in

\noindent In the hyperbolic setting, we begin with a (conjugacy class of) Fuchsian group $\Gamma$
uniformizing a point of $T_g^s$.  Specifying also a collection of horocycles, one horocycle about each puncture
of
$F=F_g^s$ (called a ``decoration'') furthermore uniquely determines
a point $\tilde\Gamma\in \widetilde{\cal T}_g^s$ by definition. We shall next describe the ``convex hull construction'', which assigns
to
$\tilde\Gamma\in \widetilde{\cal T}_g^s$ a corresponding point 
$(G_{\tilde\Gamma},w_{\tilde\Gamma})\in U_g^s$; this assignment is effectively computable as we shall see.

\vskip .1in

\noindent Here is a sketch of the convex hull construction from [6].  One may identify the
open positive light-cone $L^+$ in Minkowski three-space with the space of all horocycles in the hyperbolic plane.  Via this
identification, we find a $\Gamma$-invariant set $B$ of points in $L^+$ corresponding to the decoration, where we regard $\Gamma$ as
acting via Minkowski isometries.  One can show that $B$ is discrete in $L^+$ and consider the convex hull $H$ of $B$ in the vector
space structure underlying Minkowski space.  The extreme edges of the resulting
$\Gamma$-invariant convex body $H$ project to a collection of disjointly embedded arcs $\alpha _{\tilde\Gamma}$ connecting
punctures; furthermore, each component of $F-\cup{\alpha _{\tilde\Gamma}}$ is simply connected, and we say that
$\alpha_{\tilde\Gamma}$ ``fills'' $F$.  Given any arc family $\alpha$ filling $F$, we may define a subset
$$C(\alpha )=\{ \tilde\Gamma\in \widetilde{\cal T}_g^s:\alpha_{\tilde\Gamma}~{\rm is~homotopic~to}~\alpha \} .$$
The Poincar\'e dual of the cell decomposition $F-\cup\alpha _{\tilde\Gamma}$ of $F$ is a fatgraph $G$ embedded as a spine of $F$.
An explicit formula in terms of Minkowski geometry (which will be given in $\S$3) for the ``simplicial coordinates'', gives a
positive metric $w_{\tilde\Gamma}$ on $G$.

\vskip .1in

\noindent In the conformal setting, the effective construction maps $G_g^s\to (M_g^s\times{\bf R}_+^s)$; in contrast in the hyperbolic
setting, the effective construction (namely, the convex hull construction) maps $\tilde M_g^s\to G_g^s$ in the opposite direction!  
Just
as the conformal setting has a non-computable inverse $(M_g^s\times{\bf R}_+^s)\to G_g^s$, there is a
non-computable (or at least, very difficult to compute) inverse $G_g^s\to \tilde M_g^s$ in the hyperbolic setting.  (In fact, the paper
[7] is dedicated to the study of exactly these ``arithmetic problems'', which may be thought of as the computation of the hyperbolic
geometric data from the combinatorial data.)

\vskip .1in

\noindent The conformal and hyperbolic treatments of the cell decomposition of decorated Teichm\"uller space
are thus ``inverses'' in this sense, and each setting has its difficult theorem: surjectivity of the effective construction.
There is no known way to use one such difficult theorem to prove the other.  (One can, however, show that there is bounded distortion
from identifying simplicial coordinates with the fatgraph metric, but we shall take this up elsewhere.)

\vskip .1in

\vskip .1in

\noindent Thus, the difficult theorem in decorated
Teichm\"uller theory is that the putative cells $C(\alpha )\subseteq\widetilde{\cal T}_g^s$ are in fact cells. 
This putative cellularity is proven in [6] (independent of the Jenkins-Strebel theory) by introducing
an ``energy functional'' on a Eulidean space containing ${\widetilde{\cal T}}_g^s$ and analyzing its gradient flow in order to apply the
Poincar\'e-Hopf Theorem.  Technical details of the extension of this theorem to the bordered case are discussed in $\S$3.

\vskip .1in

\noindent There is further structure in the hyperbolic setting (for instance, global coordinates on $\widetilde{\cal
T}_g^s$ coming from Minkowski lengths, which has no analogue in the conformal setting), as is discussed in [6] and extended to
the bordered case in $\S$3.

\vskip .2in

\noindent {\bf Acknowledgements.}~It is a pleasure to thank Francis Bonahon for
helpful comments and Dennis Sullivan for many useful suggestions, insights, and corrections. 

\vskip .3in

\noindent {\bf 1. Two versions of decorated moduli space} 

\vskip .2in

\noindent Having defined in the Introduction the Teichm\"uller and moduli spaces for $r=0$, both decorated and classical,
this section is dedicated solely to the definitions of our two versions of decorated moduli spaces in the bordered
case.

\vskip .1in

\noindent  
Let us henceforth assume that
$r\neq 0$, enumerate the (smooth) boundary components of $F$ as $\partial _i$, where $i=1,\ldots ,r$, and set $\partial=\cup\{
\partial _i\} _1^r$.  
\vskip .1in

\noindent Let $Hyp(F)$ be the space of all hyperbolic metrics on
$F$ with geodesic boundary and define the first version of {\it moduli space} to be
$$M=M(F)= \bigl [Hyp(F)\times (\prod _{1}^r\partial _i)\bigr ]/PF,$$
where $PF$ denotes the equivalence relation of push-forward by orientation-preserving diffeomorphism
$$f_*\bigl (\Gamma ,(\xi _i )_1^r)\bigr ) = (f_*(\Gamma ),(f(\xi _i)_1^r),$$ where  $\Gamma\mapsto f_*(\Gamma )$ is the usual
push-forward of metric on $Hyp(F)$.

\vskip .2in

\noindent {\bf Remark} We shall not require in this paper a version of the corresponding Teichm\"uller space, but comment on it here
for completeness. In certain mathematical circles, the ``standard'' definition  of the Teichm\"uller space
$T=T(F)$ of the bordered surface is as follows:  Choose a point $\xi _i\in\partial _i$, for each $i=1,\ldots ,r$ and let   
$T(F)$  be the quotient $Hyp(F)$ by push-forward by diffeomorphisms
fixing each
$\xi _i$, where the diffeomorphisms are isotopic to the identity with the isotopy likewise required to fix each $\xi _i$.
$T$ can be shown to be homeomorphic to an open ball of real dimension $6g-6+4r+2s$.
$PMC$ acts on $T$ by push-forward with a quotient $T/PMC$ which is non-naturally homeomorphic to the moduli space $M$ just
defined.  To define a homeomorphism $T/PMC\to M$ requires using the geometry of $\Gamma$ to concoct well-defined basepoints in
the universal cover of each $\partial _i$, and it is delicate.

\vskip .2in

\noindent Letting $\ell _i(\Gamma)$ denote the hyperbolic length of $\partial _i$ for $\Gamma\in Hyp(F)$, define the first version of
{\it decorated moduli space} to be
$$\tilde M=\tilde M(F)= \bigl\{ \bigl (\Gamma ,(\xi _i )_1^r,(t_i)_1^r\bigr ):\Gamma\in Hyp(F), \xi _i\in\partial _i, 0<t_i< \ell
_i(\Gamma ), i=1,\ldots ,r\bigl\} /PF,$$
where $PF$ denotes push-forward by diffeomorphisms on $(\Gamma ,(\xi _i )_1^r )$ as before, extended by the trivial action on $(t_i
)_1^r$.  Thus, a point of $\tilde M$ is represented by $\Gamma\in Hyp(F)$ together with a pair of points $\xi _i\neq p_i$ in each
$\partial _i$, where $p_i$ is the point at hyperbolic distance $t_i$ along $\partial _i$ from $\xi _i$ in the orientation
on $\partial _i^*$ as a boundary component of $F^*$.
There is one special case, namely
$g=0=s=r-2$, so
$F$ is an annulus; in this case, we define
$\tilde M(F)$ to be the collection of all configurations of two distinct labeled points in a circle of some radius.

\vskip .1in

\noindent This completes the definition of the first version, and we turn now to the second
version.  Begin with a smooth surface $F$ with smooth boundary, choose one distinguished
point
$d_i\in\partial _i$ in each boundary component, and set
$D=\{ d_i\} _1^r$.  (We could take $d_i=\xi _i$, for instance, but it would be confusing notation in the sequel.)
Define a {\it quasi hyperbolic
metric} on
$F$ to be a hyperbolic metric on
$F^\times =F-D$ so that each $\partial _i^\times =\partial _i-\{ d_i\}$ is totally geodesic, for each $i=1,\ldots ,r$, and set $\partial
^\times=\cup\{
\partial _i^\times\} _1^r$.  To explain this, consider a hyperbolic metric on a  once-punctured annulus $A$ and the simple geodesic
arc
$a$ in it asymptotic in both directions to the puncture which separates the two boundary components; the induced metric on a
component of $A-a$ gives a model for the structure near a component of $\partial ^\times $; see Figure~2.  

\vskip .1in

\noindent The {\it decorated Teichm\"uller space} is the space
$\widetilde{\cal T}=\widetilde{\cal T}(F)$ of all quasi hyperbolic metrics on $F-D$, 
where we furthermore specify for each $d_i$ a segment of a horocycle
centered at $d_i$, modulo push-forward by diffeomorphisms of $F-D$ which are isotopic to the identity; diffeomormphisms of $F$ act
trivially on hyperbolic lengths of horocylic segments by definition. We shall see in Theorem~1 that $\widetilde{\cal T}$ is homeomorphic
to an open ball of dimension
$6g-6+5r+2s$. The second version of {\it decorated moduli
space} is
$\widetilde{\cal M}=\widetilde{\cal M}(F)=\widetilde{\cal T}(F)/PMC(F)$.

\vskip -.7in

\hskip 1.7in{{{\epsfysize2in\epsffile{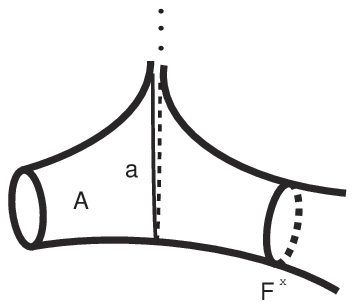}}}}

\centerline{{\bf Figure 2}~~{The model for $F^\times$.}}

\vskip .2in

\noindent Fix some quasi hyperbolic metric on $F$.  In the homotopy class of $\partial _i$ is a unique
separating geodesic $\partial _i^*\subseteq F$.  We may excise from $F-\cup\{\partial _i^*\}_1^r$ the components containing
points of $\partial ^\times$ to produce a surface $F^*$, which inherits a hyperbolic metric.  In the special
case of an annulus, the surface $F^*$ collapses to a circle.  As a point of notation, $\tilde\Gamma\in\widetilde{\cal T}$
has its underlying hyperbolic metric given by a conjugacy class of Fuchsian group $\Gamma$ for
$F^*$.

\vskip .3in

\noindent {\bf 2. Arc families}

\vskip .2in

\noindent Define an {\it (essential) arc} in $F$ to be a smooth path $a$
embedded in $F$ whose endpoints lie in $D$ and which meets $\partial F$ transversely, where
we demand that
$a$ is not isotopic rel endpoints to a path lying in $\partial F -D$.  Two arcs
are said to be {\it parallel} if there is an isotopy between them
which fixes $D$ pointwise.  
An {\it arc family} in
$F$ is the isotopy class of a collection of disjointly embedded
essential arcs in
$F$, no two of which are parallel.  

\vskip .1in
 \noindent If $\alpha$ is a collection of arcs representing an arc family in $F$, we shall say that an embedded
arc or curve
$C$ in
$F$ {\it meets
$\alpha$ efficiently} if there are no bigons in $F$ complementary to $\alpha\cup C$.

\vskip .1in

\noindent Suppose that $\alpha$ is an arc family in $F$ so that each component of
$F-\cup \alpha $ is either a polygon or a once-punctured polygon; in this case, we shall say that
$\alpha $ {\it quasi fills} the surface $F$.  In the extreme case that each component is a triangle
or a once-punctured monogon, then
$\alpha$ is called a {\it quasi triangulation}.

\vskip .3in

\noindent {\bf 3. Lambda lengths and simplicial coordinates }

\vskip .1in

\noindent  We begin with a global coordinatization of $\widetilde{\cal T}$ and recall that if
$h_0$ and
$h_1$ are two horocycles in the hyperbolic plane, then their {\it lambda length} is $\sqrt{2~{\rm exp}~ \delta}$, where
$\delta$ is the signed hyperbolic distance between $h_0$ and $h_1$ (and the sign is positive if and only if $h_0$ and $h_1$
are disjoint).  Via the canonical identification of the open positive light-cone in Minkowski space with the space of all
horocycles in the hyperbolic plane (cf. [6;$\S$1]), the square of the lambda length is simply the negative of the Minkowski inner
product (cf. [6; Lemma 2.1]).

\vskip .2in

\noindent{\bf Theorem 1} ~~\it Fix any quasi triangulation $\tau$ of $F$.  Then the assignment of lambda lengths
defines a real-analytic homeomorphism
$$\widetilde{\cal T}_{g,r}^s\to {\bf R}_+^{\tau\cup\partial ^\times}.$$\rm

\vskip .2in

\noindent{\bf Proof}~As in the proof of Theorem~3.1 of [6], we may choose a triangle complementary to $\tau$ in $F$ and a triple
of rays in the light-cone in Minkowski space and uniquely realize a triple of putative lambda lengths on a chosen triangle of
$F-\tau$ with a triple of points in these rays.  One may then uniquely and inductively construct lifts of adjacent triangles to
the light-cone realizing the putative lambda lengths in order to produce a tesselation.  Finally (and following Poincar\'e), one
explicitly constructs the underlying Fuchsian group as the group of hyperbolic symmetries of this tesselation, which leaves
invariant the corresponding set of horocycles by construction.~~~~ \hfill {\it q.e.d.}

\vskip .2in

\noindent{\bf Side-Remark}~In fact, one can give a representation of $PMC$ as a group of rational functions acting on lambda
lengths as follows.  Fix a quasi triangulation $\tau$, and adopt lambda lengths coordinates for $\widetilde {\cal T}$ with
respect to $\tau$.  If $f\in PMC$, then there is a sequence of ``elementary transformations'' (i.e., replacing one diagonal of a
quadrilateral with the other) which carries $f(\tau )$ to $\tau$.  In order to describe the action of $f$, we must calculate the
length of the other diagonal from the one.  
If $e$ is an arc in a decorated quasi triangulation $\tau$ which
separates two triangles with respective sides $a,b,e$ and $c,d,e$ and $f$ is the diagonal other than $e$ of the quadrilateral
with sides $a,b,c,d$, then it can be shown  [6; Proposition~2.6] that
$f=e^{-1} (ac+bd)$, where we have identified an arc with its lambda length for convenience. Such a transformation of lambda
lengths is called a ``Ptolemy transformation'' owing to its kinship with the classical theorem of Ptolemy on 
Euclidean quadrilaterals which inscribe in a circle.  Thus, after a permutation induced by $f$, the action of
$f$ on lambda lengths with respect to $\tau$ is given by a composition of Ptolemy transformations.  For instance, we calculate the
representation of the braid groups in the addendum to [6].  There is also a simple expression for the Weil-Petersson K\"ahler
two-form in lambda lengths [9; Theorem 3.3.6] in the case of surfaces without boundary; the invariance of this expression
under Ptolemy transformations, which devolves to a simple calculation, shows that this same expression provides a $PMC$-invariant two-form
on
$\widetilde{\cal T}$ and hence a two-form on $\widetilde{\cal M}$ itself.  Finally, [6; $\S$6] shows that ``centers'' of cells,
corresponding to setting the lambda lengths identically equal to unity, are uniformized by Fuchsian groups $\Gamma < PSL_2({\bf R})$ that
are arithmetic in the sense that there is a representative of the conjugacy class with $\Gamma <PSL_2({\bf Z})$.  An analogous
statement holds in the bordered case as well, where there is a representative $\Gamma\in PSL_2({\bf R})$ so that each
$\gamma\in\Gamma$ lies in $PSL_2({\bf Z})$ {\sl except} those hyperbolics corresponding to boundary geodesics;
there is further interesting arithmetic structure associated with these exceptional covering transformations which deserves
further study.

\vskip .2in

\noindent  For the second coordinatization, recall that if $e$ is an arc in a decorated quasi triangulation $\tau$ which
separates two triangles with respective sides $a,b,e$ and $c,d,e$, then the {\it simplicial coordinate} of $e$ is
$$E={{a^2+b^2-e^2}\over{abe}}+{{c^2+d^2-e^2}\over{cde}},$$
where we have identified an arc with its lambda length for convenience.  In the special case that $e$ bounds a once-punctured
monogon, define its simplicial coordinate to vanish; in the special case that $e\in\partial ^\times$, it bounds a
triangle on only one side, say with edges $a,b,e$, and we define its simplicial coordinate to be
$E=2~{{a^2+b^2-e^2}\over{abe}}$ (and are thus taking the usual simplicial coordinate in the double of $F$).

\vskip .1in

\noindent Fix a quasi triangulation $\tau$ of $F$, and define the subspace
$$\eqalign{
\tilde C(\tau )&=\{ (\vec{y},\vec{x})\in{\bf R}^{
\partial ^\times}\times ({\bf R}_+\cup\{ 0\} )^{\tau }:
~{\rm there~are~no~vanishing~cycles~or~arcs}\}\cr
&\subseteq {\bf R}^{\tau\cup\partial ^\times} ,\cr
}$$
where the coordinate functions are taken to be the simplicial coordinates (rather than lambda lengths as in Theorem~1).
By ``no vanishing cycles'', we mean there is no essential simple closed
curve $C\subseteq F$ meeting a representative $\tau$ efficiently so that
$$0=\sum _{p\in C\cap\cup\tau} x_p,$$ where $p\in C\cap a$, for
$a\in\tau $, contributes to the sum the coordinate $x_p$ of
$a$.  By ``no vanishing arcs'', we mean there is no essential simple arc
$A\subseteq F$ meeting $\tau$ efficiently and properly embedded in $F$
with its endpoints disjoint from $D$
so that
$$0=\sum _{p\in A\cap \partial ^\times} y_p~~+~~\sum
_{p\in A\cap\tau} x_p,$$
where again the $x_p$ and $y_p$ denote the coordinate on $a$
at an intersection point $p=A\cap a$ or $p=C\cap a$ for $a\in\tau$.  (There are
always two terms in the former sum.)
We may think of this as a convex constraint on $\vec{y}$ given
$\vec{x}$.

\vskip .1in

\noindent Fix a quasi triangulation $\tau $ of $F$, and let
$\sigma _i$ denote the triangle in $F$ complementary to
$\tau$ which contains
$\partial _i^\times$, for $i=1,\ldots ,r$.
We shall require the following analogue of Lemma~5.2 of [6].

\vskip .2in

\noindent{\bf Lemma 2}~~\it Suppose that $(\vec y,\vec x)\in \tilde C(\tau )$.  If the strict triangle inequality on the lambda
length of
$\partial _i^\times$ in $\sigma _i$ holds, for each $i=1,\ldots, r$, 
then all
three strict triangle inequalities hold on the lambda lengths of any triangle complementary to $\tau$.
\rm

\vskip .2in 

\noindent {\bf Proof}~Adopt the notation in the definition of simplicial coordinates for the lambda lengths near an edge $e$.  If
$c+d\leq e$, then $c^2+d^2-e^2\leq -2cd$, so the non-negativity of the simplicial coordinate $E$ gives
$0\leq cd  [(a-b)^2-e^2]$, and we find a second edge-triangle pair so that the triangle inequality fails.
Define an {\it arc of triangles} $(t_j)_1^n$ to be a collection of triangles complementary to $\tau$ so that
$t_j\cap t_{j+1}=e_j$, for each $j=1,\ldots ,n-1$, and likewise define  
a {\it cycle of triangles} when $t_j\cap t_{j+1}=e_j$, for each $j=1,\ldots , n$, taking the index $j$ to be cyclic (so
that
$t_{n+1}=t_1$).  In either case, if the edges of $t_j$ are $\{ e_j,e_{j-1},b_j\}$, for $j=1,\ldots ,n$, then the collection $\{
b_j\} _1^n$ is called the {\it boundary} of the cycle, and the edges $\{ e_j\}$ are called the 
{\it consecutive edges} of the cycle.
It follows that
if there is any such
triangle $t$ so that the triangle inequalities do fail for $t$, then there must be a cycle of triangles of such failures or an
arc of triangles of such failures whose boundary begins and ends with elements of $\partial ^\times$.  The former possibility is
untenable since if
$e_{j+1}\geq b_j+e_j$, for $j=1,\ldots n$, where we again identify an arc with its lambda length, then upon summing and canceling
like terms, we find $0\geq \sum _{j=1}^n b_j$, which is absurd since lambda lengths are positive.~~~~~\hfill {\it q.e.d.}

\vskip .2in

\noindent As discussed in the Introduction, Theorem 5.4 of [6] is our version of the reverse Jenkins-Strebel Theorem in decorated
Teichm\"uller theory (and is proved independently of the usual Jenkins-Strebel Theorem), and it gives a $PMC$-invariant cell
decomposition of
$\widetilde{\cal T}$. In effect, $\tilde\Gamma$ gives rise to a quasi filling arc family $\alpha _{\tilde\Gamma}$ via the convex hull
construction; fixing the topological type of $\alpha
_{\tilde\Gamma}$ and varying
$\tilde{\Gamma}$ gives a cell in the decomposition of $\widetilde{\cal T}$. 
The extension of the convex
hull construction and  its associated cell decomposition to bordered surfaces is given by  

\vskip .2in

\noindent {\bf Theorem 3} ~~\it There is a real-analytic homeomorphism of the
decorated Teichm\"uller space $\widetilde{\cal T}$ of
$F$ with $\bigl [\bigcup _\tau \tilde C(\tau )\bigr ]/\sim$, where 
$(\tau _1,\vec y_1,\vec x_1)\sim (\tau _2,\vec y_2,\vec x_2)$ if
$\vec y_1=\vec y_2$ and  
$\vec x_1$ agrees with $\vec x_2$ on 
$\tau _1-\{a\in\tau _1:x^1=0\}=\tau _2-\{a\in\tau _2:x^2=0\}$,
where $x^j$ denotes the $\vec x$ coordinate on $a$, for $j=1,2$.
Indeed, a point $\tilde\Gamma\in\widetilde{\cal T}$ gives rise to the quasi filling arc family $\alpha _{\tilde\Gamma}$ via the
convex hull construction as well as a tuple of simplicial coordinates $(\vec y,\vec x)\in\tilde C(\tau)$ for any quasi triangulation
$\tau\supseteq\alpha_{\tilde\Gamma}$, where $\vec x$ vanishes on $\tau-\alpha _{\tilde\Gamma}$.
\rm

\vskip .2in

\noindent{\bf Proof}~
The proof closely follows that of Theorem 5.4 of [6] in the double ${\cal F}$ of $F$, where one takes the
convex hull in Minkowski space of the set of all horocycles in $F$ to produce an invariant convex body.

\vskip .1in

\noindent A more technical discussion of the extension to our present situation is as follows.  As in
Theorem~5.4, the argument involves an ``energy functional'', which is defined exactly as in the proof of Theorem~5.4 (p. 322) but on the
double ${\cal F}$.  From the very definition, notice that there is no new contribution to the functional for the edges lying in $\partial
^\times$ since the coupling equations automatically hold for these edges in the double; there is therefore no need for further
computational elaboration beyond the cases considered in [6; Cases 1-8, pp. 324-327].  Furthermore, the cycle of triangles argument
in Claim~1 (p. 322) again extends to a cycle or arc of triangles argument (just as in Lemma~2).  The remaining proof of Theorem 5.4
now holds verbatim.  Consider a face of the convex body corresponding to a triangle complementary to a quasi
triangulation; if this triangle contains points of $\partial ^\times$, then it is not necessarily the case that the support plane of
this face is elliptic, but any other support plane is either elliptic as in [6; p. 320] or parabolic as in [6; p.
336].~~~~~\hfill{\it q.e.d.}

\vskip .3in

\noindent In particular, for a ``generic'' point $\tilde\Gamma\in\widetilde{\cal T}$, the arc family $\alpha _{\tilde\Gamma}$
arising from the convex hull construction is a quasi triangulation.

\vskip .1in

\noindent To close this section, we recall an elementary fact
which is useful in $\S$6, where we study certain asymptotic problems, and to formulate this fact, we recall the notion of ``h-length''.
Consider a decorated ideal triangle with lambda lengths of consecutive edges given by $a,b,e$ (where we use these same symbols also for
the edges themselves).  The geodesics $a,b$ cut off a finite horocyclic segment from the horocycle opposite $e$, and a calculation
[6; Proposition 2.8] shows that half the horocyclic length, or {\it h-length} of this segment is given by $e/ab$.  Notice that the
simplicial coordinate of an edge is by definition a six-term linear combination of the h-lengths near the edge.  

\vskip .1in

\noindent Given a cycle of triangles in an ideal triangulation of $F$, each arc in its boundary is opposite a well defined horocyclic
segment, called an ``included'' horocyclic segment.  

\vskip .2in

\noindent{\bf Lemma 4}\it ~~~Given any cycle of triangles, the sum of the simplicial coordinates of the consecutive edges is twice the
sum of the included h-lengths.\rm

\vskip .2in

\noindent
\noindent {\bf Proof}~~Adding the six-term linear relations for consecutive edges in a cycle of triangles,
the formula follows from elementary cancellation.~~~~\hfill{\it q.e.d.}

\vskip .3in

\noindent {\bf 4. Equidistant points to horocycles}

\vskip .2in

\noindent In  [7], we studied equidistant points to horocycles in the hyperbolic plane and next recall the results of
attendant elementary and explicit calculation.

\vskip .2in

\noindent{\bf Lemma 5}~~\it  Given three horocycles $h_0,h_1,h_2$ with distinct centers, let $\lambda _j$ denote
the lambda length of $h_k$ and $h_\ell$, where $\{ j,k,\ell\} =\{ 0,1,2\}$. 

\vskip .1in

\noindent {\bf a)} {\rm [7; Proposition 2.3]}~~There is a point $\zeta$ in the hyperbolic plane which is equidistant from $h_0$,
$h_1$, and
$h_2$ if and only if $\lambda _0,\lambda _1,\lambda _2$ satisfy all three possible strict triangle inequalities; in this
case, $\zeta$ is unique, and fixing the centers and varying only the decorations, all points of the hyperbolic plane arise.
Finally, the exponential $\rho$ of the common hyperbolic distance from $\zeta$ to $h_0$, $h_1$, or $h_2$ is given
by
$$\rho ^2={{2\lambda _0^2\lambda _1^2\lambda _2^2}\over{(\lambda _0+\lambda _1+\lambda _2)
(\lambda _0+\lambda _1-\lambda _2)(\lambda _0+\lambda _2-\lambda _1)(\lambda _1+\lambda _2-\lambda _0)}}.$$

\vskip .1in

\noindent {\bf b)} {\rm [7; Proposition 2.5]}~~If $\sigma$ is the geodesic 
connecting the centers of
$h_k$ and  $h_\ell$, then the signed hyperbolic length of the horocyclic segment between $\sigma\cap h_k$ and the
central projection of $\zeta$ to $h_k$ is given by
$${{\lambda _k^2+\lambda _\ell ^2-\lambda _j^2}\over{4\lambda _j\lambda _k\lambda _\ell}},$$
where the sign is positive if and only if $\sigma$ does not separate $\zeta$
from the center of $h_j$.
In particular, if $\sigma$ does separate $\zeta$ from the center of $h_j$, then $\lambda _j^2>\lambda _k^2+\lambda
_\ell^2$.

\vskip .1in

\noindent {\bf c)}~Suppose that $e$ is a diagonal of a decorated quadrilateral where the lambda lengths satisfy all three strict
triangle inequalities on each triangle complementary to $e$, and let $\zeta ,\zeta '$ denote the corresponding equidistant points from
part a). Choose an endpoint
of $e$ and centrally project
$\zeta ,\zeta '$ to the horocycle centered at this endpoint.  The simplicial coordinate of $e$ vanishes if and only if
these projections coincide.\rm

\vskip .2in

\noindent  {\bf Proof}~The reader is referred to [7] for the computational proofs of parts a) and b).  Part c) follows directly
from part b) and the definition of simplicial coordinates as in [7].~~~\hfill{\it q.e.d.}

\vskip .1in

\noindent In order to guarantee the existence of equidistant points and apply Lemma~5, we must pass to a strong
deformation retract $\widehat{\cal M}=\widehat{\cal M}(F)\subseteq
\widetilde{\cal M}$ of $\widetilde {\cal M}$.  The subspace
$\widehat{\cal M}$ is most easily defined as the $PMC$-quotient of another space 
$$\widehat{\cal T}=\widehat{\cal T}(F)=\bigl  [\bigcup _\tau \hat C(\tau )\bigr ]/\sim,$$
where $\sim$ is as in Theorem~3, and for membership in $\hat C(\tau )$ we demand not only that there are no vanishing cycles
or arcs, but we also require that for any triangle
$t\subseteq F$ complementary to
$\tau$, the lambda lengths on the edges of $t$ satisfy all three possible strict triangle inequalities.
This defines the subspace $\widehat{\cal M}\subseteq\widetilde{\cal M}$.  

\vskip .2in

\noindent {\bf Lemma 6}\it ~~$\widehat{\cal M}\subseteq \widetilde{\cal M}$ is a strong deformation retraction.\rm

\vskip .2in

\noindent {\bf Proof}~Again consider the triangle $\sigma _i$ containing $\partial _i^\times$ which is complementary to some quasi
triangulation $\tau$.  If the triangle inequality on lambda lengths fails for $\partial _i^\times$ in $\sigma _i$, then we may
simply decrease the lambda length of $\partial _i^\times$ in order to ensure that the strict triangle inequality holds on the
resulting lambda length of $\partial _i^\times$ in $\sigma _i$, for each $i=1,\ldots ,r$.  According to Lemma~2, there can then
be no failure of strict triangle inequality on the lambda lengths for any triangle complementary to $\tau$.  This homotopy of
$\widetilde{\cal T}$ to $\widehat{\cal T}$ descends to give the asserted strong deformation retraction of 
$\widetilde {\cal M}$ to $\widehat {\cal M}$.~~~~~\hfill{\it q.e.d.}

\vskip .2in

\noindent We close this section with a lemma which will be useful in $\S$6.

\vskip .2in

\noindent {\bf Lemma 7}~~\it Suppose $\alpha$ is an arc family arising from the convex hull construction for some point of
$\widehat{\cal M}$ and an arc in
$\alpha$ has corresponding lambda length $e$ and simplicial coordinate $E$.  Then $eE\leq 4$.\rm

\vskip .2in

\noindent {\bf Proof}~Observe that in the notation of the definition of simplicial coordinate, we have
$$eE=(a^2+b^2-e^2)/ab+(c^2+d^2-e^2)/cd.$$
Since the underlying decorated hyperbolic structure lies in $\widehat{\cal M}$, each triple $a,b,e$ and $c,d,e$ satisfy all three
possible triangle inequalities.  By the Euclidean law of cosines, the righthand side of the previous equation can be interpreted as
twice the sum of two cosines, and is therefore at most four.~~~~\hfill{\it q.e.d.}

\vskip .3in

\noindent {\bf 5. Isomorphism of the two versions}

\vskip .2in

\noindent Given a generic point $\tilde\Gamma\in\widehat{\cal T}$, let $\tau=\alpha _{\tilde\Gamma}$ denote the quasi triangulation
arising from the convex hull construction.  We may take $\tau$
to consist entirely of $\Gamma$-geodesics.  

\vskip .1in

\noindent Insofar as $\tau$ is a quasi triangulation, $\partial _i^\times$ lies in the
frontier of an ideal triangle (which was called $\sigma _i$ before) of $F-\tau$, and we may choose a lift $t_0$ of this ideal
triangle to the universal cover
$U$ of
$F^\times$. $U$ is a proper subset of the hyperbolic plane which is bounded by lifts of the various $\partial
_i^\times$. Consider the orbit of
$t_0$ under the primitive hyperbolic covering transformation
$\gamma$ corresponding to $\partial _i^*$ whose axis $G$ meets $t_0$, and adopt the notation illustrated in Figure~3 for the edges and
vertices in
$\{\gamma ^j(t_0)\} _{-\infty}^\infty$: the vertices of $t_0$ are $u_0,u_1,v_0$, the edge $c_0$ of $t_0$ covers $\partial
_i^\times$ with endpoints $u_0,u_1$, the remaining edges of $t_0$ are $a_0$ which has endpoints $u_0,v_0$ and
$b_0$ which has endpoints $u_1,v_0$, and
$z_j=\gamma ^j (z_0)$, for
$z=u,v,a,b,c$ and $j\in{\bf Z}$.  

\vskip .1in

\noindent Notice that at each such vertex there is a well-defined horocycle derived from the decoration.
Since $\partial _i^*$ inherits an orientation from that of $F$,
$\gamma \in\{\gamma  ^{\pm 1}\}$ can be well-defined, and we suppose that $a_1$ separates $t_0$ from the attracting fixed point
at infinity of
$\gamma $.  

\vskip -.5in
 
\hskip 1.3in{{\epsfysize3in\epsffile{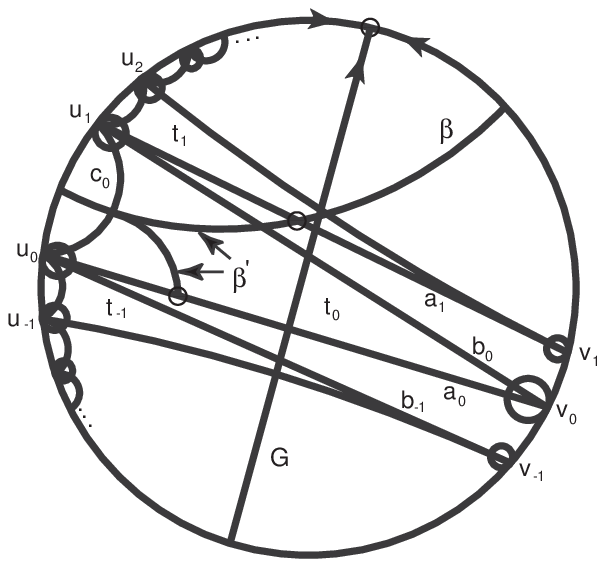}}}

\vskip .1in

\centerline{{\bf Figure 3}~~{The orbit of $t_0$.}}

\vskip .2in

\noindent The horocycles $h_0,h_1,k_0$ centered respectively at $u_0,u_1,v_0$ admit a unique equidistant point $\zeta$
according to part a) of Lemma~4 since $\tilde \Gamma\in\widehat{\cal T}$.  Since the horocycles
$h_0$ and $h_1$ both cover the same horocyclic arc, there is a curve $\beta$ in the
hyperbolic plane of possible points equidistant to them, and $\beta$ is simply the perpendicular bisector of $c_0$; that is,
$\beta\cap c_0$ is equidistant to the horocycles $h_0$ and $h_1$, and $\beta$ is perpendicular to $c_0$;  
$\beta$ is furthermore asymptotic to
$v_0'=\delta v_0$, for some $\delta\in\Gamma$.  (The arc $\beta '$ will be explained later.)

\vskip .1in

\noindent We may define a projection $$\pi :U\to \partial _i^*$$ 
as
follows: In each $t_j$, $\pi$ is induced by central projection from $v_j$, and on the component of $U-\cup\{
t_j\}_{-\infty}^\infty$ lying between
$t_{k-1}$ and $t_k$, $\pi $ is induced by central projection from $u_k$; these combine to give a continuous surjection $\pi
=\pi _\tau$ which will be useful later.

\vskip .1in
 
\noindent For any geodesic $a$ in the hyperbolic plane which is 
disjoint from the interior of $t_0$, let $H(a)$ denote the
half-plane of $a$ containing $t_0$.

\vskip .2in

\noindent {\bf Lemma 8}~~\it In the notation above, $\zeta\in H(a_1)\cap H(b_{-1})$.  \rm

\vskip .2in

\noindent {\bf Proof}~~Suppose for instance that $\zeta\notin H(a_1)$.  By part b) of Lemma 5, the triangle inequality on
squares of lambda lengths fails on the edges of $t_0$, and by part c) of Lemma~5, each of the
equidistant points of the triangles of $U-\tau$ lying in $H(a_1)$ which contain $u_1$ as vertex also must lie in the
complement of
$H(a_1)$.  There is thus a cycle of triangles of such failures, and the argument of Lemma~2
again applies to derive a contradiction.  (Thus, whereas we apply the logic of Lemma 5.2 in  [6] to the lambda lengths
themselves in Lemma~2, here we apply part of this logic to the squares of the lambda lengths.)  The argument for
$H(b_{-1})$ is analogous.~~~\hfill{\it q.e.d.}

\vskip .2in

\noindent Let $f_i$ denote the Dehn twist along $\partial _i^*$ with induced
action
$\tilde f_i$ on $\widehat{\cal T}$, for $i=1,2,\ldots ,r$.  Pushing forward by $f_i$ has the effect of moving $v_0'$ to
$\delta v_{\pm 1}$.  Thus, in the
${\bf Z}$-orbit $\tilde f_i^j(\tilde\Gamma)$ generated by this Dehn twist, there is some least $j$, call it $J$, so that
$v_0'$ lies in the complement of $H(a_1)$.  

\vskip .1in

\noindent Passing now to the $PMC$-orbit of $\tilde\Gamma$, we may replace $\tilde\Gamma$ by $\tilde f_J(\tilde \Gamma )$ as these
represent the same point of $\widehat{\cal M}$. 

\vskip .1in

\noindent The universal
cover $\tilde{\cal F}$ of the double ${\cal F}$ of $F$ is obtained by gluing together
copies of
$U$ along the lifts of the various $\partial _i^\times$ and may be identified with the hyperbolic plane.  Let $\omega$ denote the
M\"obius transformation which interchanges $u_0$ and
$u_1$ and maps $\tilde\iota (v_0)$ to $v_0$, where $\tilde\iota $ denotes reflection in $c_0$; that is,
$\tilde\iota$ is the lift of the canonical involution $\iota :{\cal F}\to {\cal F}$ which setwise fixes $c_0$.

\vskip .1in

\noindent Notice that if $\zeta\notin H(c_0)$, then $\omega (\zeta )\in H(a_0)\cap H(b_0)$ by part b) of Lemma~5 since the
simultaneous failure of two triangle inequalities (one weak and one strict inequality on the squares of lambda lengths) among positive
numbers is absurd.  Thus, $\zeta\notin H(c_0)$ implies that $\omega(\zeta )\in t_0$.  

\vskip .1in

\noindent In any case,
$\zeta '=\{
\zeta,\omega (\zeta )\}$ is a well-defined point in
$U$ which lies on the piecewise-smooth arc 
$$\beta '=[\beta\cap H(a_1)\cap H(c_0)]\cup [\omega (\beta \cap \tilde\iota (t_0))]$$ illustrated in Figure~3.
In fact, $\zeta '$ lies in the interior of $\beta '$ by Lemma 8 and the previous paragraph.

\vskip .1in

\noindent
Define
$$\eqalign{
A_i&=A_i(\tilde\Gamma )=\pi _\tau (a_1\cap\beta ')=\pi _\tau(a_0\cap\beta '),\cr
B_i&=B_i(\tilde\Gamma )=\pi _\tau (\zeta '),\cr
C_i&=C_i(\tilde\Gamma )=\pi _\tau (\beta \cap G),\cr
}$$
where $G$ is the axis of $\gamma$, i.e., $G$ is the lift of $\partial _i^*$ depicted in Figure~3.

\vskip .2in

\noindent {\bf Lemma 9}~~\it $B_i$ is a well-defined point in
$\partial _i^*-\{ A_i\}$, and $C_i$ is a well-defined point of $\partial _i^*$, for $i=1,\ldots ,r$, where
$\tau$ is any completion to an ideal triangulation of the arc family $\alpha _{\tilde\Gamma}$ associated to $\tilde\Gamma$ via
the convex hull construction.\rm 

\vskip .2in

\noindent {\bf Proof}~Suppose that the convex hull construction assigns an arc family $\alpha=\alpha_{\tilde\Gamma}$ to
$\tilde\Gamma$ which is not a quasi triangulation.  Complete $\alpha$ in any manner to a quasi triangulation $\tau$, and extend
simplicial coordinates by setting them to zero on the arcs of $\tau -\alpha$.  According to part c) of
Lemma~5, the vanishing of simplicial coordinates precisely guarantees the independence of $B_i$ from the choice of $\tau$.
$C_i$ is well-defined independent of this choice by construction.

 ~~~~~\hfill{\it
q.e.d.}

\vskip .2in

\noindent Notice that if all three horocycles centered at
the vertices of $t_0$ cover the same horocyclic arc in $F$, for instance when $r=1$, then as we vary the length of this
horocyclic arc,  the equidistant point $\zeta$ does not change.  On the other hand, when we vary the
length of the horocylic arc corresponding to the common horocycles at $u_0$ and $u_1$, the {\sl distance} to the
equidistant point varies continuously in any case.  

\vskip .1in

\noindent Define $\delta _i$ to be the signed distance from $\zeta '$ to $h_0$, where the
sign is positive if and only if $h_0$ and $h_1$ are disjoint; 
the exponential ${\rm e}^{\delta _i}$ of this distance $\delta _i$ is
expressed in lambda lengths in part a) of Lemma~5.   
Furthermore, $\delta _i$ gives an abstract coordinate on $\beta '$ in any case for
$\delta ^-_i<
\delta _i< \delta ^+_i$, where
$\delta ^+_i$ corresponds to $\zeta ^+=\beta '\cap a_1$ and $\delta ^-_i$ corresponds to $\zeta ^-=\beta '\cap a_0$.
(One can in fact compute $\delta _i^\pm$ in terms of lambda lengths and the entries of $\gamma$, but the formula is
complicated and unnecessary.)

\vskip .1in

\noindent Define the ``oriented distance'' $d_o(X,Y)$ between two distinct points $X,Y\in\partial _i^*$ to be the distance along
$\partial _i^*$ from $X$ to $Y$ in the orientation as a boundary component of $F^*$.  Notice that $d_o(X,Y)+d_o(Y,X)=\ell _i$ for any 
distinct $X,Y\in\partial _i^*$.  We shall also say that the point $X\in\partial _i^*$ is at ``signed oriented distance'' $d$ from
$Y\in\partial _i^*$, if $d_o(Y,X)=d$ or equivalently $d_o(X,Y)=\ell _i-d$.

\vskip .1in

\noindent Define $\xi _i=C_i\in\partial _i^*$; we shall define $p _i$ by specifying its signed oriented distance from $\xi _i$.
Define $\ell _i^+=d_0(A,B)$ and $\ell ^-=i-d_o(B,A)$, so $\ell _i=\ell _i^+-\ell _i^-$,
and let
$f_i: [\delta _i^-,\delta _i^+]\to [\ell _i^-,\ell _i^+]$ be the orientation-preserving affine homeomorphism.
Define $p _i\in\partial _i^*$ to be
the point in
$\partial _i^*$ at signed oriented distance from $\xi =C _i$ given by 
$$g_i(\delta _i)=({\delta ^+_i-\delta ^-_i}){{f_i(\delta _i)}}.$$

\vskip .2in

\noindent{\bf Theorem 10}~~\it The mapping $F\mapsto F^*$ together with the assignment of points $p_i,\xi _i\in\partial
_i^*$, for
$i=1,2,\ldots ,r$, gives a well-defined  real-analytic homeomorphism
$$\Psi: \widehat {\cal M}\to \tilde{ M}$$
provided $F$ is not the annulus $F_{0,2}^0$. \rm 

\vskip .2in

\noindent {\bf Proof}~~ We have already proven that the mapping $\Psi :\widehat{\cal M}\to \tilde M$ is well-defined, and it is
obviously real-analytic.  We must still show that $\Psi$ is a bijection.  

\vskip .1in

\noindent To this end, first notice that by construction, a change in the hyperbolic length of the horocyclic arc at $d_i$ leaves $\xi
_i$ invariant, moves
$p _i$ monotonically relative to $\xi _i$, and leaves invariant each $p_j,\xi _j$, for $i\neq j=1,\ldots ,r$.
Furthermore, the distance of $p _i$ from any point varies at a uniform rate
${\ell _i}={{\ell _i^+-\ell_i^-}}$ as a function of $\delta _i$.

\vskip .1in

\noindent In contrast, the
Fenchel-Nielsen deformation [13] along $\partial _i^*$ moves the endpoint of $\beta$ in $H(c_0)$ monotonically along the circle at
infinity by definition, so 
$\xi _i=C_i$ moves monotonically along $\partial _i^*$ in the deformation parameter by construction.  Thus, the Fenchel-Nielsen
deformation monotonically varies
$\xi _i$ along $\partial _i^*$ while fixing each $p_j,\xi _j$ for $j\neq i$; $p _i$ also moves under this Fenchel-Nielsen deformation,
and there is a compensatory deformation of horocyclic segment at $d_i$ so that $p_i$ and $\xi _i$ both move keeping $d_o(\xi _i,p _i)$
invariant. 

\vskip .1in 

\noindent As to surjectivity, first observe that the composition
$\widehat{\cal M}\mapright{\Psi} \tilde M \to M$ is surjective for $F\neq F_{0,2}^0$ since given a hyperbolic metric, we may always
adjoin a unique hyperbolic surface of type $F_{0,1}^{0,1}$ to get a quasi hyperbolic metric on $F$ lying in $\widehat{\cal M}$.
Since the compensated Fenchel-Nielsen deformations along $\partial _i^*$ attain all possible $\xi _i\in\partial _i^*$, it
remains only to alter the horocyclic segment at $d_i$ to vary $p _i$. 

\vskip .1in

\noindent To see that $\Psi$ is injective,
again use that there is a unique hyperbolic surface of type $F_{0,1}^{0,1}$ to conclude that if $\Psi(\tilde\Gamma
_1)=\Psi (\tilde\Gamma _2)$, then the quasi hyperbolic metrics underlying $\tilde\Gamma _1$ and $\tilde\Gamma _2$ must
agree up to a Fenchel-Nielsen deformation on $\partial _i^*$.  Again, $\xi  _i$ is 
monotone in the deformation parameter, and $d_o(\xi _i,p_i)$ is monotone in the size of the
horocycle by construction, completing the proof of injectivity.~~~~~\hfill{\it q.e.d.}

\vskip .2in

\noindent We close this section with two corollaries, which are not required in the sequel but
serve to better illuminate aspects of the homeomorphism $\Psi$.  There is the following immediate corollary to the proof of
Theorem~10.

\vskip .2in

\noindent {\bf Corollary 11}~~\it {\bf a)}~For each $i=1,\ldots ,r$, scaling the lambda length of each edge
$e\in\tau\cup\partial ^\times$ by a factor $t\in{\bf R}$ raised to the power of the number of ends of $e$ asymptotic to $d_i$
fixes $\xi _i$, moves $p _i$ uniformly around $\partial _i^*$, and leaves invariant the underlying hyperbolic
metric as well as leaving invariant each
$p_j,\xi _j$, for $i\neq j=1,\ldots ,r$.  

\vskip .1in

\noindent {\bf b)} Fix some $\hat\Gamma\in\sigma [\alpha ]\subseteq \widehat{\cal M}$ with underlying conjugacy class of Fuchsian group
$\Gamma$, and let $h_i\in{\bf R}$ denote the hyperbolic length of the horocyclic segment at $d_i$, for $i=1,\ldots ,r$.  There are
$h_i^\pm\in{\bf R}_+$ so that $h_i^-<h_i<h_i^+$.  All values in $(h_i^-,h_i^+)$ occur, and $h_i^\pm$ depend only on $\Gamma$, i.e.,
$h_i^\pm$ are independent of $h_j$, for $i\neq j=1,\ldots ,r$. \rm~~~\hfill{\it q.e.d.}

\vskip .1in

\noindent The geometrically natural procedure of moving only $p_i$ (fixing $\xi _i$ as well as the other data) can be formulated as an
action of a groupoid on $\tilde M$ as follows.  Let $(0,1)= I^\circ\subseteq I=[0,1]$ denote the unit intervals.   Define an associative
operation on $x,y\in I^\circ$ to take value $x+y$ provided $x+y<1$ and to be undefined otherwise.  This endows $I^\circ$ with
the structure of a groupoid.  In the same manner, $(I^\circ)^r$ is a groupoid with operation induced by vector sum. 
Thus, $I^\circ$ is a sub-groupoid of the additive group $S^1=I/(0\sim 1)$, and $(I^\circ )^r$ is a sub-groupoid of $(S^1)^r$.
$(I^\circ )^r$ acts as a groupoid
on $\tilde M$ in the natural way, where $\vec x=(x_i)_1^r$ maps $p_i$ to the point at oriented distance
$d_o(\xi _i,p_i)+ x_i\ell _i$ provided $x~+~d_o(\xi _i,p_i)/\ell _i<1$ for each $i=1,\ldots ,r$ (and leaves all other data unchanged),
and the action is undefined if any of these conditions fail to hold.
There is again a diagonal action of $I^\circ$ on $\tilde M$ defined in the same manner, and we shall let $\tilde M/I^\circ$ denote the
quotient.

\vskip .1in

\noindent Furthermore, let $\widehat {\cal M}/{\bf R_+}$ denote the quotient of $\widehat {\cal M}$ by the natural homothetic action of
${\bf R}_+$ on lambda lengths, or equivalently (by homogeneity of simplicial coordinates) the homothetic action on simplicial
coordinates;  an orbit of this ${\bf R}_+$-action corresponds geometrically to altering the decoration by moving each horocycle by a
common hyperbolic length fixing its center.

\vskip .2in

\noindent{\bf Corollary 12}~~\it  The homeomorphism $\Psi$ of Theorem 8 descends to a real-analytic homeomorphism $\Psi :\widehat{\cal
M}/{\bf R}_+\to
\tilde M/I^\circ$.\rm

\vskip .2in

\noindent{\bf Proof}~Using Corollary 11a, since simplicial coordinates are a homogeneous function of lambda lengths (of degree -1) by
definition and since each coordinate $\rho _i$ is likewise a homogeneous function of lambda lengths (of degree +1) by part a) of Lemma~5,
we may calculate that the speed ${d\over{d~{\rm \ell n}t}} g_i(\delta _i +\ell n~t)$ is constant equal to the hyperbolic length $\ell
_i=\ell _i^+-\ell _i^-$ of
$\partial _i^*$, for each $i=1,\ldots ,r$.~~~~\hfill{\it q.e.d.}

\vskip .3in

\noindent {\bf 6. Circle actions and the arc complex}

\vskip .2in

\noindent We first discuss circle actions on and quotients of $\tilde M$ and $M$ which will be required in the sequel.  There is a
natural ${\bf R}_+^r$-action on $\tilde M$, where $(t_i)_1^r\in {\bf R}_+^r$ replaces the vector $(\ell _i)_1^r$ of hyperbolic
lengths of boundary components by $(t_i\ell _i)_1^r$.  There is a corresponding diagonal ${\bf R}_+$-action on $\tilde M$, and we shall
let $\tilde M/{\bf R}_+$ denote the quotient; this action descends to $M$, and we shall also let $M/{\bf R}_+$ denote the quotient.  

\vskip .1in

\noindent There is a natural $(S^1)^r$-action on $\tilde M$, where the $i^{\rm th}$ factor $S^1$ moves only $p_i$ and $\xi _i$
uniformly along $\partial _i^*$ at a speed given by the hyperbolic length $\ell _i$ of $\partial _i^*$.  This action is {\sl not}
fixed-point free, and the quotient
$\tilde M/(S^1)^r$ is homotopy equivalent to the usual moduli space $M_g^{r+s}$ of the unbordered surface $F_{g,0}^{r+s}$.  On the
other hand, there is a subgroup
$(S^1)^{r-1}$ which preserves the relative positions of pairs of $\xi _i$ which evidently does act without fixed-points.  
This action descends to a well-defined action of $(S^1)^r$ on $M$ itself, where the $i^{\rm th}$ factor uniformly moves only $\xi _i$.

\vskip .1in

\vskip .2in

\noindent Now, let us inductively build a simplicial complex $Arc' (F)$, where
there is one
$p$-simplex
$\sigma(\alpha )$ for each arc family $\alpha $ in
$F$ of cardinality $p+1$.  The simplicial structure of $\sigma (\alpha )$ is the natural one, where
faces of 
$\sigma (\alpha )$ correspond to sub arc families of $\alpha $.  We begin with a vertex in $Arc'(F)$
for each isotopy class of essential arc in $F$ to define the 0-skeleton.  Having thus inductively
constructed the $(p-1)$-skeleton of $Arc'(F)$, for $p\geq	 1$, let us adjoin a $p$-simplex for each
arc family $\alpha $ consisting of $(p+1)$ essential arcs, where we identify the proper faces of
$\sigma (\alpha )$ with simplices in the $(p-1)$-skeleton in the natural way.  Identifying the open standard $p$-simplex with the
collection of all real-projective $(p+1)$-tuples of positive reals assigned to the vertices, $Arc'(F)$ is identified with the
collection of all arc families in
$F$ together with a real-projective weighting of non-negative real numbers, one such number assigned to each component of
$\alpha$.

\vskip .1in

\noindent $PMC(F)$ acts on $Arc'(F)$ in the natural way,
and we define the {\it arc complex} of $F$ to be the quotient $$Arc(F)=Arc'(F)/PMC(F).$$
If $\alpha $ is an arc
family in $F$ with corresponding simplex $\sigma (\alpha )$ in $Arc'(F)$, then we shall let $ [\alpha
]$ denote the
$PMC(F)$-orbit of $\alpha $ and $\sigma  [\alpha ]$ denote the quotient of $\sigma (\alpha )$
in $Arc(F)$. 

\vskip .1in

\noindent (The sphericity conjecture from [8] is that $Arc(F)$ is piecewise-linearly homeomorphic to the sphere of dimension
$6g-7+4r+2s$.  Furthermore, [3] introduces and studies a new topological operad whose
underlying  spaces are homeomorphic to open subspaces of $Arc(F)$.)

\vskip .2in 

\noindent For each $i=1,\ldots ,r$, there is a natural $S^1$-action on $Arc(F)$ itself as follows.  Suppose that the projective
class of a positive weight $w$ on $\alpha$ represents a point of $Arc(F)$ for some quasi-filling arc
family $\alpha$.  Imagine fattening each component arc 
$\alpha$ to a band whose width is given by the weight of $w$ assigned to the component; thus, a positively weighted arc family is
regarded as a collection of bands of positive widths running between the boundary components, and conversely.
Certain of the arcs in $\alpha$ meet some $\partial _i$ at the point $d_i$ (or rather, are asymptotic to $d_i$ in the
surface $F^\times$ with totally geodesic boundary); let $a_1,\ldots ,a_k$ denote this set of arcs meeting $d_i$ with corresponding
weights
$w_1,\ldots ,w_k$.    Now, given $t\in S^1=I/(0\sim 1)$, alter the family of bands corresponding to $\{ a_i\} _1^k$ in a neighborhood of
$\partial _i^\times$ as follows:
twist a total width $t\sum _1^k w_i$ of the bands to the right around the boundary as illustrated in Figure~4.  The resulting bands
then determine a weighted arc family, and this gives
$Arc(F)$ the structure of an $(S^1)^r$-space.

\hskip .2in{{{\epsfysize2.5in\epsffile{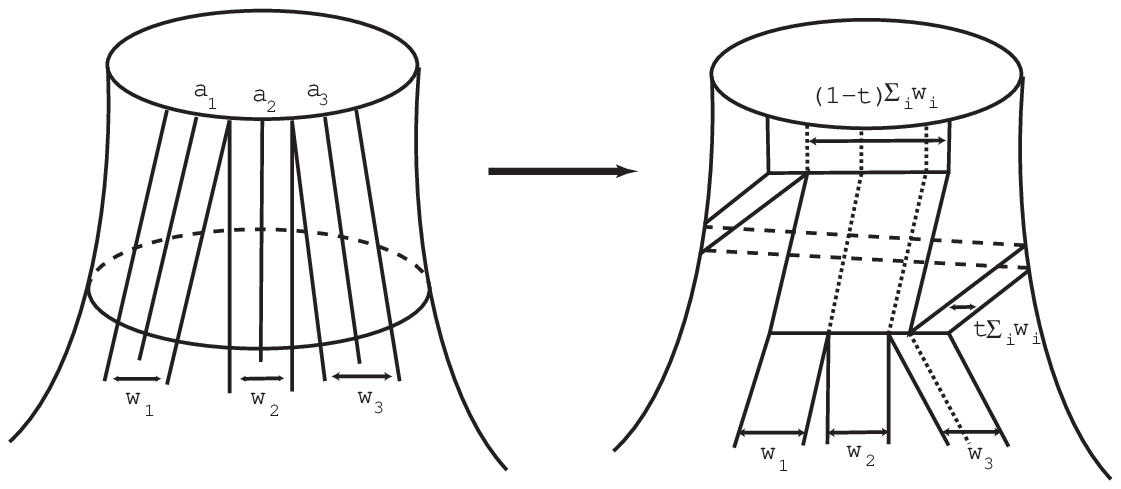}}}}

\vskip .1in

\centerline{{\bf Figure 4}~~{Circle action on $Arc$.}}

\vskip .2in

\noindent A subspace of $Arc(F)$ of special interest to us here and elsewhere in the general theory is
$$Arc_\#(F)=\{ \sigma [\alpha ]: \alpha ~{\rm quasi-fills}~F\} .$$
The $(S^1)^r$-action on $Arc(F)$ preserves $Arc_\#(F)$.
It is useful to have a ``deprojectivized'' $Arc_\#(F)$, and we define $dArc_\# (F)=Arc_\#(F)\times {\bf R}_+$
which we identify with all positive-real weightings on the components of an arc family which quasi fills $F$.

\vskip .1in

\noindent For each quasi triangulation $\tau$, define
$$\tilde{D}(\tau)=\{ \vec x:~{\rm there~are~no~vanishing~cycles~for}~\vec x\} ,$$
and notice that by definition 
$$dArc_\#(F)=\biggl (\bigl [\bigcup _\tau \tilde D(\tau )\bigr ]/\sim \biggr )/PMC(F),$$
where the equivalence relation is as in Theorem~3.  

\vskip .2in

\noindent {\bf Lemma 13}~~\it The natural mapping $q:\widehat{\cal M}\to dArc_\# (F)$ induced by
projection $(\vec y,\vec x)\mapsto \vec x$ is a homotopy equivalence.  The fibers of $q$ are
homeomorphic to the interior of a cone over a finite-sided convex polyhedron.
\rm

\vskip .2in

\noindent {\bf Proof}~Let $\vec 0$ denote the $r$-dimensional vector consisting entirely of entries zero, so if
$\vec x\in\tilde D(\tau )$, then $(\vec 0,\vec x)\in\tilde C(\tau )$.  As noted before Lemma~2, $q^{-1}(\vec x)$ is convex and
hence strong deformation retracts to $(\vec 0,\vec x)$.~~~~~\hfill{\it q.e.d.}

\vskip .2in

\noindent {\bf Remark}~Forgetting the simplicial coordinates on $\partial ^\times$ is a violent
operation: Fixing the simplicial coordinates $\vec x$ on each arc in some quasi triangulation $\tau $ and altering only the
simplicial coordinates $\vec y$ of arcs in $\partial ^\times$ changes the underlying hyperbolic metric in an extremely
complicated and non-computable way (cf. [7]).  One characteristic shared by $[\vec y_1,\vec x],[\vec y_2,\vec x]\in \tilde C(\tau
)$ is as follows: For any essential simple closed curve $C$ in $F$, we may assume that $C$ meets $\tau$ efficiently and consider the
sum over $z\in C\cap \tau$ (with multiplicity and without sign) of the simplicial coordinate of $\vec x$ on the arcs meeting $C$ at
$z$, which is a kind of length of $C$ (as seen from the horocycles) as in [7; Lemma~1.2]. For each curve $C$, these lengths coincide
for
$[\vec y_1,\vec x]$ and
$[\vec y_2,\vec x]$.  A similar combinatorics in [1] captures the hyperbolic lengths of geodesics. 

\vskip .2in

\vskip .2in

\noindent {\bf Theorem 14}~\it For any bordered surface $F\neq F_{0,(2)}^0$, $Arc_\# (F)$ is proper homotopy equivalent
to $M(F)/{\bf R}_+$ as $(S^1)^r$-spaces.\rm

\vskip .2in

\noindent{\bf Proof}~~According to Lemma~13, the fiber of $q: \widehat {\cal M}\to dArc_\#$ is the interior of a cone over a
finite-sided convex polyhedron, and the fiber of the induced map $\widehat{\cal M}/{\bf R}_+\to Arc_\#$ is an open finite-sided convex
polyhedron.  In particular, $Arc_\#, dArc_\#, \widehat{\cal M}, \widehat{\cal M}/{\bf R}_+$ all have the same homotopy type.
Likewise, $M,M/{\bf R}_+,\tilde M, \tilde M/{\bf R}_+$ all have the same homotopy type.  Using the homeomorphism $\Psi :\widehat{\cal
M}\to \tilde M$ of Theorem~10, it follows that
 $Arc_\# (F)$ is indeed homotopy equivalent to $M(F)/{\bf R}_+$.  Notice that these two spaces furthermore have the same dimension.

\vskip .1in

\noindent To explicitly describe the map $Arc_\#(F)\to M/{\bf R}_+$, given a projective weight $w$ on an arc family $\alpha$
representing a point of $\sigma [\alpha ]$, regard $w$ as projective simplicial coordinates on $\alpha\cup\partial ^\times$, where the
simplicial coordinate on each component of $\partial ^\times$ is taken to be zero.  This point of $\widehat{\cal M}$ gives rise via the
construction in the proof of Theorem~8 to a point of $\tilde M$; we finally forget the decoration and projectivize the hyperbolic
lengths of the boundary geodesics to describe a point of $M/{\bf R}_+$.  The standard Fenchel-Nielsen deformations about the boundary
curves act on each of the spaces $\widehat{\cal M},\tilde M, M/{\bf R}_+$; taking the point $\xi _i$ as the initial point of an arc
family meeting $\partial _i^*$ for each $i$ identifies the corresponding $(S^1)^r$-action on $Arc_\#$ with the
action described before on $Arc(F)$.  Thus, the homotopy equivalence between $Arc_\#$ and $M/{\bf R}_+$ is indeed a map of
$(S^1)^r$-spaces.

\vskip .2in

\noindent It remains only to prove properness.  To this end, since we mod out by the homothetic ${\bf R}_+$-action on
$Arc_\# =Arc_\#(F)$, we may and shall assume that all lambda lengths are bounded below, or in other words by homogeneity, that all
simplicial coordinates are bounded above.    In the same manner, since we mod out by the ${\bf R}_+$-action on $M=M(F)$, we may and
shall assume that all hyperbolic lengths of geodesic boundaries are likewise bounded above.  

\vskip .1in

\noindent To prove properness, first imagine a path tending to infinity in
$M/{\bf R_+}$.  As is well-known, there must either be 1) an essential and non-boundary parallel
curve $g$ in $F$ whose hyperbolic length is tending to zero; or 2) the hyperbolic length of some geodesic boundary component is
tending to zero.  In the former case, the lambda lengths of all arcs meeting $g$ must tend to infinity.  It follows from Lemma 7,
that the simplicial coordinates must tend to zero; that is, there is a vanishing cycle, so the corresponding path also tends to infinity
in $Arc_\#$.  

\vskip .1in

\noindent In the latter case, re-consider Figure~3.  Not drawn in the figure are the ideal polygons complementary to
$\cup \{ t_i\}_{-\infty}^{+\infty}$.  There is a unique such polygon with vertex $u_i$, for each $i$, and these polygons are generically
triangulated in the special manner where each edge has $u_i$ among its two endpoints.  To fix a particular lift, consider the
triangulated polygon $P$ with vertex $u_1$.  The hyperbolic trnasformation $\gamma$ discussed before with axis $G$ can be described as
follows: it is the composition $\gamma =\sigma\circ \tau$ of two parabolic transformations, where $\sigma$ is the parabolic fixing
$v_0$ mapping $a_0$ to $b_0$, and $\tau$ is the parabolic with fixed point $u_1$ which maps $b_0$ to $a_1$.   
Indeed, $\sigma\circ\tau$ maps $u_{-1}\mapsto u_1$ and $v_0\mapsto v_1$ by definition.  Since the lambda lengths on the edges of $t_0$
agree with those of $t_1$ (and lambda lengths are M\"obius-invariant), it follows that also $u_1\mapsto u_2$; this uniquely determines
$\sigma\circ\tau$, and is the unique M\"obius transformation $\gamma$ preserving the triangulation and mapping
$t_0$ to $t_1$.

\vskip .1in

\noindent The number of sides of $P$
is uniformly bounded in terms of only the topology of $F$; we shall refer to this fact as ``bounded combinatorics''.  Thus, the trace of
$\gamma$ is bounded away from
$\pm 2$ provided also the
$h$-lengths are bounded away from zero.  In other words, if the hyperbolic length of a geodesic boundary component tends to zero, then
some h-lengths must tend to zero, i.e., the corresponding path again tends to infinity in $Arc_\#$.

\vskip .1in

\noindent In the other direction, consider a path tending to infinity in $Arc_\#$, so there is some cycle $C$ of triangles
with vanishing simplicial coordinates.  There are again two basic cases depending upon whether the cycle $c$ in $F$ dual to $C$
satisfies: 1) $c$ is not boundary parallel; or 2) $c$ is boundary parallel.  In the latter case if $c$ is parallel to $\partial
_i$, all the simplicial coordinates along $C$ vanish, so by Lemma 4, all of the h-lengths incident on $d_i$ likewise vanish. 
Again by bounded combinatorics, the trace of $\gamma$ tends to $\pm 2$, so the corresponding path also tends to infinity in $M/{\bf
R}_+$.  In the former case, the trace of a matrix representing $c$ again has trace tending to $\pm 2$, so in either case,
the corresponding path also tends to infinity in $M/{\bf R}_+$.~~~~~\hfill{\it q.e.d.}

\vskip .2in

\noindent Together with the putative sphericity theorem, we would obtain:

\vskip .2in

\noindent {\bf Corollary 15}~\it ~~~ For any bordered surface $F\neq F_{0,(2)}^0$, the arc complex $Arc(F)$ is
a spherical compactification of $M(F)/{\bf R}_+$.\rm

\vskip .2in

\noindent {\bf Remark}~In fact, the embedding $(h_i^-,h_i^+)\to \widehat{\cal M}$ in
Corollary~11b can furthermore be shown to be proper.  As in Theorem~14, the induced homeomorphism in Corollary 12 can then 
be used to describe a compactification of $\tilde M/{\bf R}_+/I^\circ$ by a ``decorated'' arc complex
$Arc(F)\times ((S^1)^r/S^1)$, which evidently supports a natural $(S^1\times S^1)^r$-action that one can prove extends
the $(S^1\times I^0)^r$ groupoid action discussed in $\S$5.  Furthermore,
this decorated arc complex is closely related to the operads computed in [3], which naturally extend the arc
operad.  This remark will be taken up elsewhere.

\vfill\eject

\noindent {\bf Bibliography}

\vskip .2in

\noindent [1]~~ L. Chekhov and V. Fock, ``Obvervables in 3D gravity and geodesic algebras'', {\it Czech. Jour. Phys.} {\bf 49}
(1999).

\vskip .1in

\noindent [2] J. L. Harer, ``The virtual cohomological dimension of the mapping class group of an orientable surface, {\it Inv. Math.}
{\bf 84} (1986), 157-176.

\vskip .1in

\noindent [3]~~Ralph L. Kaufmann, Muriel Livernet, R. C. Penner, ``Arc operads and arc algebras'', 
preprint math.GT/0209132 (2002)

\vskip .1in

\noindent [4]~~S. Kojima, ``Polyhedral decompositions of hyperbolic manifolds with boundary'', {\it On the geometric structure of
manifolds}, ed. Dong Pyo Chi, (1990), 35-57.

\vskip .1in

\noindent [5] M. Kontsevich, ``Intersection theory on the moduli space of curves and the matrix Airy function'', {\it Comm. Math. Phys.}
{\bf 147} (1992), 1-23.

\vskip .1in

\noindent [6]~~ R. C. Penner, ``The decorated Teichm\"uller space of  punctured surfaces", 
{\it Comm. Math.  Phys.}  {\bf 113}   (1987),  299-339.

\vskip .1in

\noindent [7]~~---, ``An arithmetic problem in surface geometry'', {\it
The Moduli Space of Curves}, Birk-h\"auser (1995), eds. R. Dijgraaf, 
C. Faber, G. van der Geer, 427-466.

\vskip .1in

\noindent  [8]~~---,``The simplicial compactification of Riemann's moduli space'', 
Proceedings of the 37th Taniguchi Symposium, World Scientific (1996), 237-252.

\vskip .1in

\noindent [9]~~---, ``Weil-Petersson volumes'', {\it Jour. Diff. Geom.} {\bf 35} (1992), 559-608.

\vskip .1in

\noindent [10] ---, ``Perturbative series and the moduli space of Riemann surfaces",   
{\it Jour. Diff. Geom.} {\bf  27}  (1988),  35-53.

\vskip .1in

\noindent [11]~J. Hubbard and H. Masur, ``Quadratic differentials and foliations'', {\it Acta Math.} {\bf 142} (1979), 221-274.

\vskip .1in

\noindent [12] K. Strebel, {Quadratic Differentials}, {\it Ergebnisse der Math.} {\bf 3:5}, Springer-Verlag, Heidelberg (1984).

\vskip .1in

\noindent [13]~S. Wolpert, ``On the symplectic geometry of deformations of a hyperbolic surface'', {\it Ann. Math.}
{\bf 117} (1983), 207-234.

\vskip .1in

\bye